\documentclass[a4paper, 12pt]{amsart}
\usepackage{amssymb, amsmath, amsthm, color, MnSymbol}
\usepackage[colorlinks=true, breaklinks=true, linkcolor=black, citecolor=blue, urlcolor=red]{hyperref}
\usepackage[francais,english]{babel}
\usepackage{graphicx}

\usepackage{tikz-cd}
\usepackage{tikz}
\usetikzlibrary{arrows.meta, shapes, positioning}
\usetikzlibrary{decorations.markings}

\pgfdeclarelayer{background}
\pgfdeclarelayer{nodelayer}
\pgfdeclarelayer{edgelayer}
\pgfsetlayers{background,edgelayer,nodelayer,main}
\tikzset{
  none/.style={}
}
\tikzset{every loop/.style={min distance=12mm, looseness=4}} 

\theoremstyle{plain}
\newtheorem{theorem}{Theorem}[section]
\newtheorem{proposition}[theorem]{Proposition}

\newtheorem{theoremIntro}{Theorem}

\theoremstyle{definition}

\theoremstyle{remark}

\newtheorem{example}[theorem]{Example}

\newcommand{\cA}{\mathcal A}
\newcommand{\cAL}{\mathcal A_\Lambda}

\newcommand{\C}{\mathbb{C}}
\newcommand{\CKA}{Cuntz-Krieger algebra}

\DeclareMathOperator{\Mod}{Mod}

\newcommand{\N}{\mathbb{N}}
\newcommand{\cO}{\mathcal O}
\newcommand{\cOL}{\mathcal O_\Lambda}
\newcommand{\onto}{\twoheadrightarrow}
\newcommand{\cP}{\mathcal P}

\DeclareMathOperator{\Rep}{Rep}

\newcommand{\Sz}{Szyma\'{n}ski}
\newcommand{\tbf}{\textbf}
\newcommand{\ti}{\tilde}
\newcommand{\Z}{\mathbb Z}

\title{A new perspective on spectra of quantum spaces}

\author{Arnaud Brothier}
\date{}
\address{Arnaud Brothier\\	University of Trieste, Department of Mathematics, via Valerio 12/1, 34127, Trieste, Italy and School of Mathematics and Statistics, University of New South Wales, Sydney NSW 2052, Australia}
\email{arnaud.brothier@gmail.com\endgraf
		\url{https://sites.google.com/site/arnaudbrothier/}}
\begin{document}
\maketitle

\begin{abstract}
We consider a class of C$^*$-algebras $C(X_q)$ associated with quantum spaces such as spheres, projective spaces, and lens spaces. 
We introduce a non-self-adjoint operator algebra $A$ together with an explicit functor from the category of representations of $A$ to that of $C(X_q)$.
We then construct explicitly a family of one-dimensional representations of $A$ that parametrise the entire spectrum of $C(X_q)$.
\end{abstract}

\section*{Introduction}
Vaughan Jones found unexpectedly an efficient way to construct unitary representations of certain groups while searching for reconstruction of conformal field theories \cite{Jones17}, see the survey \cite{Brothier20}.
It was first applied to the groups $F,T,V$ of Richard Thompson \cite{Cannon-Floyd-Parry96}.
This method was then adapted to the Cuntz(--Dixmier) C$^*$-algebras $\cO_n$: any representation of $\cO_n$ can be constructed using an explicit inductive limit and a representation (often of finite dimension) of a much simpler algebra 
 \cite{Dixmier64, Cuntz77, Brothier-Jones19}.
Powerful method were then developed to classify a large part of the spectrum of $\cO_n$ and in particular in obtaining a sequence of moduli spaces in \cite{Brothier-Wijesena25,Brothier-Wijesena26,Brothier-Wijesena24a}.
These methods were then substancially adapted to the algebraic context of the Leavitt path algebras \cite{Brothier-Wijesena24b} and later on adapted to Cuntz-Krieger algebras (in fact to higher rank graph C$^*$-algebras) \cite{Brothier-Sims-Wijesena26}.

This article aims to present and apply these techniques to C$^*$-algebras of quantum spaces. 
These are presented noncommutative C$^*$-algebras of type $I$ obtained by deforming algebras of function on topological spaces $X$. 
They are interpreted as spaces of functions over a fictive ``noncommutative topological space'' $X_q$ and thus often denoted $C(X_q)$.
We will consider Vaksman-Soibelman quantum spheres (of odd and even dimension) and some of their orbifold such as quantum projective spaces and quantum lens spaces appearing in \cite{Vaksman-Soibelman90, Hong-Szymanski03,Brzezinski-Szymanski18,Gotfredsen-Zegers24}.

The study of C$^*$-algebras of such quantum spaces is a rich and vibrant research topic mixing operator algebraic techniques with geometric and algebraic ones. 
They first appeared in the pioneered work of Woronowicz and Podle\'{s} \cite{Woronowicz80, Woronowicz87, Podles87}. They are closely related to quantum groups as developed by Faddeev, Drinfeld, Jimbo and many other, see \cite{Drinfeld86,Jimbo85} and the survey \cite{Faddeev06}.
Deep work have been pursued in describing and studying their underlying topological and geometrical properties; for instance in computing their K-theory, K-homology Fredholm modules, Chern characters and even equipping certain of them with a Connes spectral triple structure, to cite a few \cite{Dabrowski-Landi-Masuda01, Hawkins-Landi04, Andrea-Dabrowski-Landi08,Connes94}.

Remarkably, Hong and \Sz\ have constructed graphs $L=L(X_q)$ and explicit isomorphisms $\varphi:\cO_L\to C(X_q)$ from the Cuntz-Krieger algebra $\cO_L$ of the graph $L$ to $C(X_q)$ for a number of quantum spaces $X_q$ including quantum spheres, quantum projective spaces and certain quantum lens spaces \cite{Hong-Szymanski02,Hong-Szymanski03}.
Both C$^*$-algebras $\cO_L$ and $C(X_q)$ have explicit presentations and the map $\varphi$ is described by basic formula involving the generators of $\cO_L$ and $C(X_q)$.
These powerful techniques were then extended to a larger class of lens spaces in a work of Brzezi\'{n}ski and \Sz\ \cite{Brzezinski-Szymanski18}. 
As later observed in \cite{Gotfredsen-Zegers24}, the original definition of the graph $L$ in that work required a subtle correction. 

In this present paper we consider a non-self-adjoint algebra $\cA_L$ associated to a graph $L$ which is roughly obtained by taking the adjoints of the usual generators of $\cO_L$ subject to all relations of $\cO_L$ except the one requiring them to be partial isometries.
For generic graphs $L$ it is in practice much easier to construct $\cA_L$-modules than representations of $\cO_L$ (we use the terms \emph{module} for actions of non-self-adjoint algebras and \emph{representations} for actions of C$^*$-algebras).
In particular, $\cA_L$ usually admits many finite dimensional modules while often $\cO_L$ has only infinite dimensional representations.
We define an explicit functor $\Pi:\Mod(\cA_L)\to \Rep(\cO_L)$ from $\cA_L$-modules to representations of $\cO_L$ that is essentially surjective (i.e.~for all $K\in \Rep(\cO_L)$ there is $H\in \Mod(\cA_L)$ so that $\Pi(H)\simeq K$). 
Moreover, $\Pi$ sends finite dimensional irreducible modules to irreducible representations and is fully faithful on those.
For each quantum object $X_q$ with associated graph $L=L(X_q)$, we construct explicitly a set of one-dimensional $\cA_L$-modules whose image under the functor $\Pi$ form a set of representatives for the spectrum of $\cO_L$.

\begin{theoremIntro}\label{theo:main}
Fix $n\geq 1$ and let $C(S_q)$ be the C$^*$-algebra of the quantum sphere of dimension $(2n-1)$. Identify $C(S_q)$ with the \CKA\ $\cO_L$ associated to the graph $L=L_{2n-1}$ of Hong and \Sz\ and consider our non-self-ajoint algebra $\cA_L$.
Then for each $(j,z)$ with $1\leq j\leq n$ and $z$ in the unit circle there is a one-dimensional module $M_{j,z}$ of $\cA_L$ so that the set of $\Pi(M_{j,z})$ parametrises the spectrum of $C(S_{q})$.
Similar results hold for quantum projective spaces, quantum lens spaces and for $2n$-dimensional quantum spheres this latter having its spectrum parametrised by $n$ circles and two additional isolated points.
\end{theoremIntro}

Parametrising these spectra is not new.
Indeed, Vaksman and Soibelman have given explicit representatives of the spectra of all quantum spheres \cite{Vaksman-Soibelman90}. 
Then Hawkins and Landi gave explicit representatives of the spectra of certain quantum spaces that turned out to be isomorphic to Vaksman-Soibelman quantum spheres \cite{Hawkins-Landi04}.
Another method goes as follows: a general theorem computes the primitive spectrum of any Cuntz--Krieger algebra \cite{Huef-Raeburn97}; see also \cite{Hong04} and \cite{Bates-Hong-Raeburn-Szymanski02}. 
Using remarkable isomorphisms between \CKA s and quantum spaces we deduce the primitive spectrum of $C(X_q)$ for $X_q$ being a quantum sphere, projective space or lens space.
Finally, using that $C(X_q)$ is a C$^*$-algebra of type $I$ one obtains the spectrum since this latter is isomorphic to the primitive spectrum via the obvious map $\pi \mapsto \ker(\pi)$.

The interest of our construction lies in a different point of view, based on the introduction of a non-self-adjoint algebra $\cA_L$ and the lifting $\Pi:\Mod(\cA_L)\to\Rep(\cO_L)$ of all $\cA_L$-module to representations of $\cO_L\simeq C(X_q)$.
In particular, we directly work on representations rather than on ideals and our techniques are explicit.
Somewhat surprisingly, every irreducible representation of $C(X_q)$ come from \emph{one-dimensional} modules of $\cA_L$. 
From this perspective, $C(X_q)$ has a very basic representation theory, akin to that of a \emph{commutative} C$^*$-algebra.
This situation contrasts sharply with that of non-type~I C$^*$-algebras such as the Cuntz algebra $\mathcal{O}_2$ in two generators.
It admits a non-self-adjoint operator algebra $\cA_2$.
The primitive spectrum of $\mathcal{O}_2$ consists of a single point, since $\mathcal{O}_2$ is simple, whereas the part of its spectrum corresponding to $d$-dimensional irreducible $\cA_2$-modules is parametrised by a smooth manifold of dimension $2d^2+1$ for all natural number $d\geq 1$ \cite{Brothier-Wijesena24a}.

\textbf{Plan of the article.}
In Section \ref{sec:preliminaries} we define Cuntz-Krieger algebras $\cO_L$ associated to finite directed graphs $L$.
We define our non-self-adjoint algebra $\cA_L$ and explain how each $\cA_L$-module $H$ dilates into a representation $\Pi(H)$ of $\cO_L$.
We recall some key properties of the process $H\mapsto \Pi(H)$.
In Section \ref{sec:quantum-objects} we consider one by one each C$^*$-algebra $C(X_q)$ of quantum spaces studying here.
We define a graph $L:=L(X_q)$ so that $\cO_L\to C(X_q)$.
Then we describe our non-self-adjoint algebra $\cA_L$ for this specific graph $L$ and construct explicitly a set of irreducible $\cA_L$-modules so that their image by $\Pi$ forms a set of representatives for the spectra of $\cO_L$ (and thus of $C(X_q)$).
We provide great details on the first examples that are the odd-dimensional quantum spheres.
Then we explain how to derive similar results for quantum projective spaces and quantum lens spaces that are all orbifolds of odd-dimensional quantum spheres.
Lastly, we describe the spectrum of the even-dimensional quantum spheres.

\section*{Acknowledgements}

We are grateful to Giovanni Landi for various discussions and encouragement related to this project, and for making specific comments on a previous version of this paper.

\section{Lifting representations}\label{sec:preliminaries}

If $L$ is a finite directed graph, then we define a non-self-adjoint operator algebra $\cA_L$ and the Cuntz--Krieger algebra $\cO_L$ associated to $L$.
Then we explain that given a representation $H$ of $\cA_L$ (we will say an $\cA_L$-module) we can construct a representation $\Pi(H)$ of $\cO_L$. We then provide key properties of the process $H\mapsto \Pi(H)$.
All this approach has been done for the larger class of row-finite and locally convex higher rank graphs $L$ in \cite{Brothier-Sims-Wijesena26}.
We refer the reader to this latter article for proofs and additional details.

\subsection{Finite directed graphs}
Consider a nonempty finite directed graph $L$ with vertex set $L^0$, edge set $L^1$ and the usual source and range maps $s,r:L^1\to L^0$.
An edge with same source and range is called a \emph{loop}. 
A path of length $n\geq 1$ is written $(\lambda_m,\cdots,\lambda_1)$ or $\lambda_m\cdots \lambda_1$ with $\lambda_i$ an edge so that $s(\lambda_{j+1})=\lambda(e_j)$ (hence we take the category convention: paths are written from right to left).
We set $L^n$ to be the set of all paths of length $n$ and set $L^*=\cup_{n\geq 0}L^n$ all paths including vertices that are considered as paths of length $0$.
We extend $r,s$ on paths so that $r(v)=s(v)=v$ if $v\in L^0$ and $r(\lambda_n\cdots \lambda_1)=r(e_n), s(\lambda_n\cdots \lambda_1)=s(\lambda_1)$.
A path $\lambda$ with nonzero length having same source and target is called a \emph{cycle}.
We equip $L^*$ with the usual partially defined composition and write $\lambda L^*$ for all paths of the form $\lambda\mu$ with $\mu\in L^*$. In particular, $vL^*$ denotes all paths with range $v$.
We similarly define $\lambda \Lambda^n$ and $v\Lambda^n$ for $n\geq 1$.

\begin{center}\textbf{From now all graphs are assumed to be nonempty finite and directed.}\end{center}

\subsection{Cuntz--Krieger algebras}
Let $L$ be a graph (hence directed finite and nonempty).
The \emph{\CKA}\ $\cO_L$ of $L$ is the universal unital C$^*$-algebra having for generators $x_v, v\in L^0$ and $x_\lambda, \lambda\in L^1$ satisfying:
\begin{enumerate}
\item the $x_v$ with $v\in L^0$ are mutually orthogonal projections so that $\sum_{v\in L^0} x_v =1$;
\item for all $w\in L^0$ we have $\sum_{\lambda\in w\Lambda^1} x_\lambda x_\lambda^*=x_w$;
\item $x_\lambda^*x_\lambda = x_{s(\lambda)}$ for all $\lambda\in L^1$.
\end{enumerate}
We will often write $x_\lambda$ for $x_{\lambda_n}\cdots x_{\lambda_1}$ when $\lambda=\lambda_n\cdots\lambda_1\in L^n$.
Note that some authors take the opposite convention than ours such as Hong and \Sz\ in \cite{Hong-Szymanski02}.
The Cuntz-Krieger algebra of a graph $L$ in our convention corresponds to the Cuntz-Krieger algebra of the opposite graph in the convention of Hong and \Sz.
In particular, our conventions imply that $x_\lambda$ for $\lambda\in L^1$ is a partial isometry satisfying $x_\lambda=x_{r(\lambda)} x_\lambda x_{s(\lambda)}$.

Being universal means that the representations of the C$^*$-algebra $\cO_L$ are in correspondence with the tuples $(K,X_v,X_\lambda:\ v\in L^0,\lambda\in L^1)$ so that $K$ is a Hilbert space, $X_v,X_\lambda\in B(K)$ are bounded linear operators on $K$ satisfying the 3 axioms of above.
Note that all representation $K$ of $\cO_L$ decomposes into the (orthogonal) direct sum of the $K_v:=X_v(H)$ for $v\in L^0$.
Then $X_\lambda$ can be interpreted as a map $K_{s(\lambda)}\to K_{r(\lambda)}$.

\subsection{Pythagorean algebras}

Consider a graph $L$.
Define the \emph{Pythagorean algebra} $\cP_L$ of $L$ to be the universal unital C$^*$-algebra having for generators $a_v$ with $v\in L^0$ and $a_\lambda, \lambda\in L^1$ satisfying 
\begin{enumerate}
\item the $a_v$ with $v\in L^0$ are mutually orthogonal projections so that $\sum_{v\in L^0} a_v =1$;
\item for all $w\in L^0$ we have $\sum_{\lambda\in wL^1} a_\lambda^*a_\lambda=a_w$.
\end{enumerate}
By universality, the following formula defines a surjective morphism of C$^*$-algebras:
$$\cP_L\onto \cO_L, \ a_v\mapsto x_v,\ a_\lambda\mapsto x_\lambda^*.$$
Hence, our choice of conventions force us to map $a_\lambda$ to the \emph{adjoint} of $x_\lambda$. Moreover, observe that $\cP_L$ admits the presentation of $\cO_L$ to which we have removed the last axiom (3).
To define $a_\lambda$ for a path $\lambda=\lambda_n\cdots\lambda_1$ we must reverse the order and set
$$a_{\lambda_n\cdots\lambda_1}:=a_{\lambda_1}\cdots a_{\lambda_n}.$$
We then have the contravariant formula for compositions:
$$a_{\lambda\mu}=a_{\mu} a_{\lambda} \text{ for } \lambda,\mu\in L^*.$$
Additionally, if $H$ is a representation of $\cP_L$ with operators $A_\lambda,\lambda\in L^*$, then $H$ decomposes into $\oplus_{v\in L^0}H_v$ where $H_v:=A_v(H)$ and $A_{\lambda}$ defines a bounded linear map $H_{r(\lambda)}\to H_{s(\lambda)}$.
By this we mean that the domain of $A_\lambda$ is contained in $H_{r(\lambda)}$ and its range contained in $H_{s(\lambda)}.$

This algebra was first defined by Jones for the graph $L$ with one vertex and $2$ loops in \cite{Brothier-Jones19}.
In this example $\cP_L$ is presented by two operators $a,b$ and one relation $a^*a+b^*b=1$.
The relation is reminiscent of the Pythagorean equality and thus gave the name to this algebra.

\subsection{A non-self-adjoint algebras}
Define now $\cA_L$ to be the norm-closed non-self-adjoint subalgbra of $\cP_L$ generated by the $a_v$ and the $a_\lambda$ with $v\in L^0,\lambda\in L^1$.
It is important to stress that $\cA_L$ is not self-adjoint. This changes drastically its representation theory.

To prevent confusions we call \emph{module} a continuous unital algebra morphism $\cA_L\to B(H)$ and we will call representation a continuous unital involutive algebra morphism $\cP_L\to B(H)$ or $\cO_L\to B(H)$.
Modules over $\cA_L$ and representations over $\cP_L,\cO_L$ define categories $\Mod(\cA_L)$ and $\Rep(\cP_L),\Rep(\cO_L)$ where the arrows are the usual intertwiners (continuous linear maps intertwining the algebra actions).
Note that $\Mod(\cA_L)$ and $\Rep(\cP_L)$ have same objects (each representation of $\cP_L$ restricts into an $\cA_L$-module and conversely each $\cA_L$-module extends into a representations of $\cP_L$) but $\Mod(\cA_L)$ has more arrows than $\Rep(\cP_L)$ (since an intertwiner of $\cA_L$-modules is not required to intertwine the adjoints $a_\lambda^*$).
A module will often be denoted $H$ (the carrier Hilbert space) and $A_v,A_\lambda$ the operators corresponding to $a_v,a_\lambda\in\cA_L$ where $v\in L^0,\lambda\in L^1.$
We call $a_\lambda$ or $A_\lambda$ an \emph{edge-operator}.
Similarly, we may write $K$ to express a representation of $\cO_L$ with operators denoted $X_v,X_\lambda$.

\subsection{Lifting representations}\label{sec:lift}
We now construct explicitly the so-called \emph{Pythagorean functor}
$$\Pi:\Mod(\cA_L)\to \Rep(\cO_L).$$

Fix a module $H=\oplus_{v\in L^0} H_v$. 
Consider all the pairs $(\lambda,\xi)$ with $\lambda\in L^*$ and  $\xi\in H_{r(v)}$.
Take all formal linear combination of those.
Now, mod out this vector space by the equivalence $\sim$ generated by:
$$(\lambda,\xi)\sim \sum_{\mu\in s(\lambda)\Lambda^1} (\lambda\mu,A_\mu\xi)$$
where $s(\lambda)\Lambda^1$ stands for all edges with range equal to the source of $\lambda$.
We set $[\lambda,\xi]$ to be the class of $(\lambda,\xi)$ with respect to $\sim$.
We obtain a vector space $K^0$ equal to all finite linear combinations of the classes $[\lambda,\xi]$.
There exists a unique inner product $\langle\cdot,\cdot\rangle$ on $K^0$ satisfying the following axioms:
\begin{enumerate}
\item $\langle [\lambda,\xi],[\mu,\eta]\rangle =0$ if $\lambda \not \in \mu L^*$ and $\mu\not\in \lambda L^*$;
\item $\langle[\gamma,\xi'], [\gamma\nu,\eta']\rangle =\langle A_\nu \xi',\eta'\rangle_H,$
\end{enumerate}
where $\lambda,\mu,\gamma,\nu$ are paths and $\xi,\eta,\xi',\eta'$ are vectors of $H$ in the right vertex-components.
Denote by $\Pi(H)$ or $K$ the Hilbert space completion of $K^0$ equipped with this inner product.
Set $K_v$ to be the closed linear span of the $[\lambda,\xi]$ with $\lambda\in vL^*$ and more generally set $K_\lambda$ to be the closed linear span of the $[\lambda\mu,\eta]$.
We obtain the following decompositions into orthogonal components
$$K=\oplus_{v\in L^0}K_v$$
and for each $m\geq 1$ we have
$$K_v=\oplus_{\lambda\in v\Lambda^{\leq m}}K_\lambda$$
where $v\Lambda^{\leq m}$ denotes all paths $\lambda$ with range $v$ of length $m$ or of length smaller than $m$ but that cannot be extended (i.e.~the source of $\lambda$ does not have any incoming edges).

We now define a representation $x_\lambda\mapsto X_\lambda\in B(K)$ of $\cO_L$ on $K$.
We set 
$$X_\lambda\cdot [\mu,\xi]=\begin{cases}
[\lambda\mu,\xi] \text{ if } s(\lambda)=r(\mu)\\
0 \text{ otherwise}
\end{cases}.$$
One must prove that this formula does not depend on the representative $(\mu,\xi)$ of the class $[\mu,\xi]$ and that the $X_\lambda$ do indeed satisfy the relations of the Cuntz--Krieger algebra $\cO_L$. 
The former is an elementary computation while the latter follows from the observations that $X_v$ is the orthogonal projection onto $K_v$ and $X_\lambda$ is a partial isometry with domain $K_{s(\lambda)}$ and codomain $K_\lambda$.
We emphasise the following useful formula concerning the adjoint of $X_\lambda$:
$$X_\lambda^*\cdot [\lambda\gamma,\eta]=[\gamma,\eta] \text{ and } X_\lambda^*[v,\xi]=[s(\lambda),A_\lambda\xi]$$
for $\lambda\in vL^1$ and $\gamma\in L^*$.

Consider now $\theta:H\to H'$ a morphism of $\cA_L$-modules, i.e.
$$\theta\circ A_\lambda=A'_\lambda\circ \theta \text{ for all path } \lambda$$
using obvious notations for the actions.
The formula 
$$[\mu,\eta]\mapsto [\mu,\theta(\eta)]$$
is well-defined and provide an intertwinner 
$$\Pi(\theta):\Pi(H)\to \Pi(H')$$
of the associated representations $\Pi(H)$ and $\Pi(H')$ of $\cO_L$.
Hence, $\Pi$ defines a functor from $\cA_L$-modules to representations of $\cO_L$.
We recall some important properties of $\Pi$. They can all be found in \cite[Proposition 5.1]{Brothier-Sims-Wijesena26}, see also \cite{Brothier-Wijesena24a} for similar statements on the Cuntz algebra $\cO_2$.

\begin{theorem}\label{theo:Pi}
Consider the functor $\Pi:\Mod(\cA_L)\to \Rep(\cO_L)$ constructed above.
The following assertions are true.
\begin{enumerate}
\item The functor $\Pi$ is essentially surjective, i.e.~for all representation $K$ of $\cO_L$ there exists a module $H$ satisfying $K\simeq \Pi(H)$ (where $\simeq$ means that there exists a unitary transformation that intertwines the actions of $\cO_L$);
\item If $H_1,H_2$ are module, then $\Pi(H_1\oplus H_2)$ is naturally isomorphic to $\Pi(H_1)\oplus \Pi(H_2)$. 
In particular, if $\Pi(H)$ is irreducible, then $H$ must be indecomposable (i.e.~$H$ does not contain two mutually orthogonal submodules). Moreover, if $H'\subset H$ is a (nonzero) submodule, then $\Pi(H')\simeq \Pi(H)$;
\item If $H$ is finite dimensional and indecomposable (e.g.~irreducible), then $\Pi(H)$ is irreducible;
\item If $H,H'$ are finite dimensional irreducible modules, then $H\simeq H'$ if and only if $\Pi(H)\simeq \Pi(H')$.
\end{enumerate}
\end{theorem}
To prove (1) take a representation of $\cO_L$, treat it as an $\cA_L$-module using the universal mapping $\cAL\to\cOL$  and apply $\Pi$ to it. The representation of $\cO_L$ obtained is unitary conjugate to the original one.
(2) is elementary. (3) and (4) are nontrivial statements. They were complete surprises when first found! They require a compactness argument. Hence, in general they are no longer true in the infinite dimensional case. 
We will below examples of infinite dimensional indecomposable modules $H$ so that $\Pi(H)$ is reducible.

\section{Quantum spaces}\label{sec:quantum-objects}

We recall the presentation of the universal C$^*$-algebra $C(X_q)$ of odd-dimensional quantum spheres.
Then we define the graph $L=L_{2n-1}$ of Hong and \Sz\ satisfying $C(X_q)\simeq \cO_L$.
From there we use our novel technology to construct a set of representatives of all irreducible classes of representations of $\cO_L.$ 
This involves defining the non-self-adjoint algebra $\cA_L$ and using the functor $\Pi:\Mod(\cA_L)\to\Rep(\cO_L)$.
This is done in great details. 
We explain how to easily adapt this analysis to certain orbifolds of $C(X_q)$ including quantum projective spaces and quantum lens spaces.
Then we apply a similar analysis to even-dimensional quantum spheres.

\subsection{Odd-dimensional quantum spheres}\label{sec:odd-dim-sphere}

Fix $n\geq 1$, define $X_q$ to be the quantum sphere of dimension $2n-1$ and write $C(X_q)$ for its C$^*$-algebra.
This C$^*$-algebra is a deformation of the usual commutative C$^*$-algebra of complex valued continuous functions on the real $(2n-1)$-dimensional sphere. 
The deformation is obtained via a real paremeter $q\in (0,1)$. The algebra $C(X_q)$ does not depend on $q$ (while the quantum object $X_q$ does in general, i.e.~the algebra of coefficient that is not norm-closed).
The C$^*$-algebra is isomorphic to the universal C$^*$-algebra with generators $z_1,\cdots,z_n$ subject to the relations
\begin{enumerate}
\item $z_jz_i=qz_iz_j$ if $i<j$;
\item $z_j^*z_i=qz_iz_j$ if $i\neq j$;
\item $z_i^*z_i = z_iz_i^* + (1-q^2) \sum_{j>i} z_j z_j^*$ for $1\leq i\leq n$;
\item $\sum_{k=1}^n z_kz_k^*=1.$
\end{enumerate}
Let $L=L_{2n-1}$ be the graph with edge set $\{1,\cdots,n\}$ and edge set all pairs $(j,i)$ (often denoted $ji$) with $1\leq i\leq j\leq n$.
The edge $ji$ has source $i$ and range $j$ that we express as an arrow $ji:j\leftarrow i$ or $ji:i\to j$.
Hong and \Sz\ have given an explicit isomorphism $\varphi:\cO\to C(S_q)$ defined on the generators where $\cO=\cO_L$ denotes the \CKA\ of $L$ \cite{Hong-Szymanski02}.

Our non-self-adjoint algebra $\cA$ associated to $L$ has for generators $n$ mutually orthogonal projections $a_1,\cdots,a_n$ and some edge operators $a_g$ with $g=(ji)$ for $1\leq i\leq j\leq n$ satisfying
\begin{enumerate}
\item $\sum_{i=1}^n a_i=1$;
\item $a_{g}=a_{s(g)} a_{g} a_{r(g)}$ for $g\in L^1$;
\item $\sum_{i=1}^j a_{ji}^*a_{ji}=a_j$ for all $1\leq j\leq n$.
\end{enumerate} 
We have an algebra morphism given by $a_i\mapsto x_i$ for $1\leq i\leq n$ and $a_g\mapsto x_g^*$ for $g\in  L^1$.
An $\cA$-module reduces to the data of:
\begin{itemize}
\item a Hilbert space $H$ that decomposes into $\oplus_{i=1}^nH_i$;
\item a family of operators $A_{ji}:H_j\to H_i$ satisfying that:
\item $\oplus_{i=1}^j A_{ji} :H_j\to \oplus_{i=1}^j H_j$ is an isometry.
\end{itemize}

We will often write $(H,A_g)$ to denote a module where the index $g$ is intended to vary of the edge set $L^1$.
Note that $H_i=a_i\cdot H$ and observe that if $g$ is an edge, then $A_g$ goes from $H_{s(g)}$ to $H_{r(g)}$ rather than the opposite.
If we want to define $A_\lambda$ for a path $\lambda$ the we obtain the following contravariant formula: $A_{\lambda\mu}=A_\mu\circ A_\lambda$ for composable paths $\lambda,\mu$.

\begin{example}
Assume that $n=3$. The graph $L=L_5$ is equal to:
$$
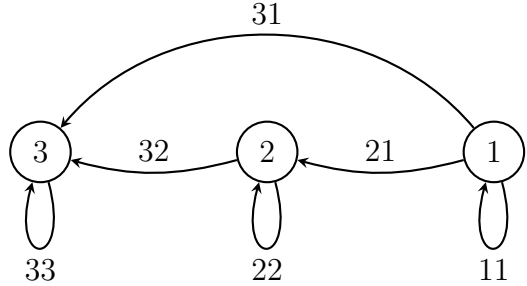
\begin{tikzpicture}[>=stealth,thick]
  \node[circle,draw,minimum size=8mm] (v3) at (0,0) {$3$};
  \node[circle,draw,minimum size=8mm] (v2) at (3,0) {$2$};
  \node[circle,draw,minimum size=8mm] (v1) at (6,0) {$1$};

  \draw[->] (v1) edge[loop below] node {$11$} (v1);
  \draw[->] (v2) edge[loop below] node {$22$} (v2);
  \draw[->] (v3) edge[loop below] node {$33$} (v3);

  \draw[->] (v1) to[bend left=15] node[above] {$21$} (v2);

  \draw[->] (v2) to[bend left=15] node[above] {$32$} (v3);

\draw[->] (v1) to[out=130,in=50] node[above] {$31$} (v3);
\end{tikzpicture}.$$
The algebra $\cA$ has for generators the orthogonal projections $a_1,a_2,a_3$ and the edge-operators $a_{11},a_{21},a_{22},a_{31},a_{32},a_{33}$ satisfying that 
\begin{itemize}
\item $a_{ji}=a_i\circ a_{ji}\circ a_j$ for all $1\leq i\leq j\leq 3$;
\item $a_1+a_2+a_3=1$;
\item $a_{11}^*a_{11}=a_1$;
\item $a_{21}^*a_{21}+a_{22}^*a_{22}=a_2$ and
\item $a_{31}^*a_{31}+a_{32}^*a_{32}+a_{33}^*a_{33}=a_3$.
\end{itemize}

An $\cA$-module is described by a Hilbert space $H$ that decomposes as 
$$H=H_1\oplus H_2\oplus H_3.$$
Then we have six edge-operators $A_{ji}:H_j\to H_i$ for $1\leq i\leq j\leq 3$ satisfying that 
\begin{itemize}
\item $A_{11}:H_1\to H_1$ is an isometry;
\item $A_{21}\oplus A_{22}:H_2\to H_1\oplus H_2$ is an isometry;
\item $A_{31}\oplus A_{32}\oplus A_{33}:H_3\to H_1\oplus H_2\oplus H_3$ is an isometry.
\end{itemize}
A submodule is then a Hilbert subspace $H'$ of $H$ that decomposes into $\oplus_i H_i'$ with $H_i':=H'\cap H_i$ and satisfying that $A_{ji}(H_j')\subset H_i'$ for any edge-operator $A_{ji}.$
\end{example}
Consider a finite-dimensional and irreducible module $(H,A_g)$ of $\cA=\cA_L$ for the graph $L=L_{2n-1}.$
Observe that $H_1$ is always a submodule and thus $H=H_1$ or $H_1=\{0\}$ by irreducibility.
Assume that $H_1=\{0\}$. 
We then obtain that $H_2\subset H$ is a submodule and thus either $H=H_2$ or $H_2=\{0\}$.
We obtain inductively that $H=H_j$ for a certain $1\leq j\leq n$.
This implies that $A_g=0$ unless $g=(jj)$: the loop at $j$. In that case $A_{jj}$ is an isometry and thus is a unitary since $H$ has finite dimension.
Moreover, all vector subspace of $H$ that is preserved by $A_{jj}$ defines a submodule.
Since $A_{jj}$ is a unitary, it is diagonalisable and therefore stabilises a line.
Irreducibility implies that $H$ must be this line and thus $\dim(H)=1$.
Therefore, $A_{jj}$ acts by multiplying by a $z\in S_1$.
Denote by $M_{j,z}$ this module: where $H=H_j=\C$, $A_{jj}=z$ and $A_g=0$ for all edge $g\neq jj$.
We have proven that all finite-dimensional irreducible module is conjugate to a certain $M_{j,z}$.
Now, it is not difficult to observe that they are pairwise inequivalent.
We deduce that 
$$\{M_{j,z}:\ (j,z)\in \{1,\cdots,n\}\times S_1\}$$
is a set of representative of the unitary classes of irreducible modules of finite dimension.
Recall that $\Pi:\Rep(\cA)\to\Rep(\cO)$ is the functor defined in Section \ref{sec:lift}.
Theorem \ref{theo:Pi} implies that the $\Pi(M_{j,z})$ parametrise a piece of the spectrum of $\cO$. 
More precisely, if $\Pi(H)$ is an irreducible representation of $\cO$ and that $\dim(H)<\infty$, then there exists a unique $(j,z)$ so that $\Pi(H)\simeq \Pi(M_{j,z}).$

Let us show that there are no other points in the spectrum of $\cO$.
Assume that $K$ is an irreducible representation of $\cO$ so that $\Pi(M_{j,z})\not\simeq K$ for all $(j,z)$.
By essential surjectivity of $\Pi$ there exists an $\cA$-module $H$ so that $\Pi(H)\simeq K$. 
Moreover, this module $H$ must be indecomposable since $\Pi(H)$ is irreducible.
Additionally, the analysis of above implies $\dim(H)=\infty$.
Assume that $H_1\neq \{0\}$. 
We must have that $\Pi(H_1)\simeq \Pi(H)$ (since $\Pi(H_1)$ is a nontrivial subrepresentation of $\Pi(H)$) and thus may assume that $H=H_1$.
We deduce that $A_{11}:H_1\to H_1$  is an isometry and all other edge operators $A_g, g\neq 11$ are equal to zero. In particular, any invariant subspace of $H_1$ for $A_{11}$ is a submodule.
The Wold decomposition gives that $A_{11}:H\to H$ decomposes, up to conjugacy, into a direct sum of the unilateral shift operator $S:\ell^2(\N)\to \ell^2(\N),\delta_n\mapsto \delta_{n+1}$ and a unitary operator $U$.
Since $\dim(H)=\infty$ and $H$ is indecomposable we must have that $A_{11}$ is conjugate to $S$.
Therefore, if $H_1\neq \{0\}$, then up to conjugacy $H=H_1=\ell^2(\N)$, $A_{11}=S$, and all other edge operator $A_g$ are equal to $0$.
Call this module $M_1$.
Assume that $H_1=\{0\}$.
In that case we can reapply a similar argument to $H_2$.
If $H_2\neq \{0\}$, then we must have, up to conjugacy, that $H=H_2=\ell^2(\N)$, $A_{22}:H_2\to H_2$ is the unilateral shift operator and all other operators $A_{g}$ are equal to zero.
Call this module $M_2$. Similarly define $M_j$ for $1\leq j\leq n$.
Continuing inductively this argument yields that our indecomposable and infinite dimensional module $H$ is necessarily conjugate to $M_j$ for a certain $1\leq j\leq n$.
We will now show that $\Pi(M_j)$ is in fact \emph{reducible}. 

We do it for $j=1$ (the other cases can be similarly obtained).
Hence, take the module $M_1=(H,A_g)$ as described above.
Set $K=\Pi(H)$ with the action of $\cO$ denoted $x_\mu\mapsto X_\mu$ as defined in Section \ref{sec:lift}.
The space $K$ is equal to the closed linear span of the $[\lambda,\xi]$ with $\lambda\in L^*$ any path and $\xi\in H_{s(\lambda)}$.
Consider the subspace $K_1$ equal to the closed linear span of the $[(11)^N,\xi]$ with $N\geq 0$ and $\xi\in H_1$ (where here $(11)^N$ denotes the path obtained by concatenating $N$ times the loop $(11)$).
Denote by $(\delta_n:\ n\in\N)$ the standard orthonormal basis of $\ell^2(\N)$ and recall that $H_1=\ell^2(\N)$.
Note that 
$$X_{11}[(11)^N,\delta_n]=[(11)^{N+1},\delta_n]$$
by definition of the action.
Moreover, $$[(11)^N,\delta_n]=\sum_{\mu \in s((11)^N)\ti L^1} [(11)^{N}\mu,A_\mu\delta_{n}]=[(11)^N(11), A_{11}\delta_n]=[(11)^{N+1},\delta_{n+1}]$$ by definition of the Hilbert space $\Pi(H)$, the fact that $(11)$ is the unique edge with range $1$ and that $A_{11}=S$.
From there we easily deduce that $X_{11}$ preserves $K_1$ and moreover that the formula
$$[(11)^N,\delta_n]\mapsto \varepsilon_{n-N}$$
defines a unitary transformation
$$U:K_1\to \ell^2(\Z)$$
where $(\varepsilon_k:\ k\in\Z)$ denotes the standard orthonormal basis of $\ell^2(\Z)$ and where $U$ conjugates $A_{11}$ to the shift operator $T:\varepsilon_k\mapsto \varepsilon_{k+1}$ of $\ell^2(\Z)$.
The shift acting on $\ell^2(\Z)$ (rather than the unilateral shift $S$ acting on $\ell^2(\N)$) is a \emph{unitary} with spectrum being the full circle.
Decompose $\ell^2(\Z)\simeq K_1$ into $V\oplus W$ where $V$ corresponds to the upper part of the circle and $W$ to the lower part of the circle.
Take now $K_V$ and $K_W$ the subrepresentations of $K$ generated by $V$ and $W$. 
We claim that $K_V$ and $K_W$ are mutually orthogonal.

Indeed, by standard results the Cuntz--Krieger algebra $\cO$ is the closed linear span of the operators $x_\lambda x_\mu^*$ with $\lambda,\mu$ paths, see \cite[Corollary 1.16]{Raeburn05} for instance.
Hence, it is sufficient to show that if $\xi\in V,\eta\in W$ and $\lambda,\mu,\alpha,\beta$ are paths, then $X_\lambda X_\mu^* \xi$ is orthogonal to $X_\alpha X_\beta^* \eta$.
Observe now that $X_\mu^* K_1=\{0\}$ if $\mu$ is not a power of the loop $(11)$.
Moreover, $X_{11}^*V=V$ and $X_{11}^*W=W.$
Hence, we may assume that $\mu$ and $\beta$ are trivial.
By definition $X_\lambda[\gamma,\zeta]=0$ if $s(\lambda)\neq r(\gamma)$.
Therefore, $X_\lambda \xi=0$ if $s(\lambda)\neq 1$ and similarly $X_\alpha\eta=0$ if $s(\alpha)\neq 1$.
Hence, we may assume that $\lambda$ and $\alpha$ have source $1$.
Note that $X_\lambda$ and $X_\alpha$ have mutually orthogonal ranges unless one is a final subpath of the other.
Hence, assume that $\lambda=\alpha\nu$ for some path $\nu$.
The path $\nu$ has source equal to the source of $\lambda$ that is $1$ and has range equal to the source of $\alpha$ that is $1$. 
Therefore, $\nu$ must be a power $(11)^N$ of the loop $(11)$.
Note that $X_\alpha^*X_\alpha$ is the othogonal projection onto the closed linear span of $[\gamma,\zeta]$ with $\gamma$ a path with range equal to the source $s(\alpha)=1$ and $\zeta\in H_{s(\alpha)}$.
In particular, $V,W$ are in the range of the projection $X_\alpha^*X_\alpha$.
We deduce that 
\begin{equation}\label{eq:inner-product}\langle X_{\alpha\nu} \xi , X_\alpha \eta\rangle = \langle X_\alpha\circ X_\nu \xi, X_\alpha \eta\rangle=\langle X_\nu \xi, \eta\rangle=\langle X_{(11)}^N\xi,\eta\rangle.\end{equation}
Since $V$ is closed under the action of $X_{(11)}$ and since $V$ is orthogonal to $W$ we deduce that \eqref{eq:inner-product} is equal to $0$.
We have proven that $K_V$ and $K_W$ are two mutually orthogonal subrepresentations of $K$.
Since $V\subset K_V,W\subset K_W$ and $V\neq \{0\}\neq W$, we deduce that $K$ is reducible.
If we slightly push this analysis we may observe that the representation $\Pi(M_1)$ can be desintegrated into the measurable field of $\Pi(M_{1,z})$ over $z\in S_1$ for the usual Lebesgue measure on the circle. Hence, $\Pi(M_1)$ is very far from being irreducible.
We have proven that $\cO$ does not admit any irreducible representation that cannot be built from a finite dimensional $\cA$-module.

\begin{theorem}\label{theo:odd-sphere}
Let $X_q$ be the quantum sphere of dimension $2n-1$ and let $C(X_q)$ be its C$^*$-algebra that we identify with $\cO_L$.
The set 
$$\{\Pi(M_{j,z}):\ 1\leq j\leq n, z\in S_1\}$$
is a set of representatives for the spectrum of $C(X_q)$.
\end{theorem}

The proof of above permits to obtain the following description of the spectrum for a specific class of \CKA s.

\begin{proposition}\label{prop:spectrum}
Consider a finite graph $L$ with $n$ vertices satisfying that 
\begin{enumerate}
\item at each vertex $v$ there is exactly one loop $\ell_v$ based at $v$;
\item there are no cycles in $L$ except product of loops.
\end{enumerate}
We write $\cA,\cO$ for the associated non-self-adjoint algebra and \CKA, respectively, and denote by $\Pi:\Mod(\cA)\to\Rep(\cO)$ the Pythagorean functor.
For each vertex $v\in L^0$ and each $z\in S_1$ define the $\cA$-module $M_{v,z}$ having for carrier Hilbert space $H=H_v=\C$, edge operator $A_g=z$ if $g=\ell_g$ and $0$ otherwise.
Then $(v,z)\mapsto\Pi(M_{v,z})$ provides an explicit parametrisation of the spectrum of $\cO$ into $n$ circles.
\end{proposition}

We will now consider other quantum spaces obtained as orbifolds of odd-dimensional quantum spheres for various finite group actions: namely quantum projective spaces and quantum lens spaces.
Their structure as Cuntz-Krieger algebras will all satisfy the assumption of Proposition \ref{prop:spectrum}.

\subsection{Orbifolds of odd-dimensional quantum spheres}

\subsubsection{Odd-dimensional quantum real projective spaces}
Fix $n\geq 1$ and $q\in (0,1)$ and let $P_q$ be the quantum real projective space of dimension $2n-1$ and write $C(P_q)$ for the corresponding C$^*$-algebra.
If $S_q$ is the $q$-deformation of the real sphere of dimension $2n-1$, then $C(P_q)$ can be realised as the fixed point subalgebra of $C(S_q)$ for the involutive transformation $x_v\mapsto -x_v$ if $v$ is a vertex and $x_\lambda\mapsto x_\lambda$ for all other paths defined on the generators of $C(S_q)\simeq \cO_L$.

Consider the graph $L_{2n-1}$ previously defined (with vertex set $\{1,\cdots,n\}$ and edge set $\{(ji):\ 1\leq i\leq j\leq n\}$ with $r(ji)=0$ and $s(ji)=i$).
Recall that the Cuntz-Krieger algebra of this graph is isomorphic to $C(S_q)$.
A procedure on the graph $L_{2n-1}$ can be performed in order to obtain that $C(P_q)$ is a graph C$^*$-algebra as well.
Indeed, define $L_{2n-1}^{(2)}$ to be the graph with same vertex set $\{1,\cdots,n\}$ than $L_{2n-1}$ but so that the edges of $L_{2n-1}^{(2)}$ are all the paths of length two of $L_{2n-1}$.
We then have that $C(P_q)$ is isomorphic to the graph C$^*$-algebra of $L_{2n-1}^{(2)}$.
Now, observe that $L:=L_{2n-1}^{(2)}$ satisfies the assumptions of Proposition \ref{prop:spectrum} and has $n$ vertices.
We deduce that the spectrum of $C(P_q)$ is in bijection with $n$ copies of circles corresponding to one-dimensional modules of $\cA_L$.

\subsubsection{Quantum lens spaces satisfying a primness condition}
We investigate quantum lens spaces in the sense of Hong and \Sz\ obtained by quantising free actions of cyclic groups on odd-dimensional spheres  \cite{Hong-Szymanski03}.
We will relax the freeness condition in the next section.
Take a natural number $n\geq 1$, a real deformation parameter $q\in (0,1)$, set $S_q$ to be the quantum sphere of dimension $2n-1$ and let $C(S_q)$ be its C$^*$-algebra.
Denote by $z_1,\cdots,z_n$ its usual generators as defined earlier.
Fix a natural number $p\geq 2$ (non-necessarily prime) and a family of natural numbers $\tbf{m}=(m_1,\cdots,m_n)$ {\bf that are all relatively prime to $p$.}
If $\theta$ is a primitive $p$th root of the unity, then the formula
$$z_i\mapsto \theta^{m_i}z_i \text{ for } 1\leq i\leq n$$
defines (by universality) an action $\Z_p\curvearrowright C(S_q)$ where $\Z_p:=\Z/p\Z$.
The fixed point C$^*$-algebra for this action is denoted $C(L_q(p;\tbf m))$ and interpreted as the C$^*$-algebra of continuous function over the quantum lens space $L_q(p;\tbf m)$. 
We will now define a graph $L''$ whose \CKA\ is isomorphic to $C(L_q(p;\tbf m))$.
First, recall that the graph $L_{2n-1}$ has vertex set $\{1,\cdots,n\}$ and edge set all pairs $(j,i)$ (often denoted $ji$) with $j\geq i$ so that the source and range of $ji$ are $i$ and $j$, respectively.
We have seen that the \CKA\ of $L_{2n-1}$ is isomorphic to the quantum sphere of dimension $2n-1$.
We will perform two transformations on $L=L_{2n-1}$ in order to obtain the desirable graph $L''$.
Define the graph $L'$ whose vertex set is $L^0\times \Z_p=\{1,\cdots,n\}\times \Z_p$ (so $np$ vertices) and edge set being $L^1\times \Z_p$: the set of pairs $(ji,m)$ with $1\leq i\leq j\leq n$ and $m\in \Z_p$ so that the source of $(ji,m)$ is $(i,m-m_i)$ and the range $(j,m)$, that is:
$$(ji,m): (i,m-m_i)\to (j,m).$$
We now define the graph $L''$ obtained from $L'$ as follows.
The vertex set of $L''$ is $\{1,\cdots,n\}$ while the edges of $L''$ from $i$ to $j$ are all the paths $\lambda$ of $L'$ satisfying
\begin{itemize}
\item the first edge of $\lambda$ is of the form $(a i,m)$ with $m=m_i$;
\item either there is only only edge or the last edge is of the form $(jb,m)$ with $m=0$;
\item all other edges are of the form $(cd,m)$ with $m\neq 0$; 
\item if $(fe,k)$ and $(hg,r)$ are different edges of the path, then $(f,k)\neq (h,r)$.
\end{itemize}
Hong and \Sz\ proved that $C(L_q(p;\tbf m))$ is isomorphic to the \CKA\ of $L''$ \cite{Hong-Szymanski03}.
Observe the following: if $i,j$ are vertices of $L''$, then there is at least one edge from $i$ to $j$ if and only if $i\leq j$.
In particular, there are no cycles except product of loops.
For each vertex $i$ there is at least one loop given by the edge $(ii,m_i)$.
It can proven that no other loops exist, see for instance \cite[Section 2]{Brzezinski-Szymanski18}.
Therefore, $L''$ has $n$ vertices and satisfy the assumption of Proposition \ref{prop:spectrum}. 

\subsubsection{General quantum lens spaces}

Recall that the C$^*$-algebra of quantum lens space $L_q(p;\tbf m)$ with $n\geq 1,p\geq 1, \tbf{m}=(m_1,\cdots,m_n)\in\N^n$ is obtained as a fixed point subalgebra of $C(S_q)$ for an action of $\Z_p$ where $S_q$ is a $(2n-1)$-dimensional quantum sphere.
We had the assumptions that the $m_i$ are all relatively prime to $p$.
This construction was generalised where now no relatively prime condition is required.
It was proven that $C^*(L_q(p;\tbf m))$ was isomorphic to a Cuntz-Krieger algebra $\cO_G$ with $G$ a precise graph (see\cite{Brzezinski-Szymanski18} and \cite{Gotfredsen-Zegers24} where the construction of the graph $G$ had to be corrected).
We will not properly define $G$ but will recall some of its keys properties.
The graph $G$ has vertex set 
$$\{v_i^s:\ 1\leq i\leq n, s\in E_i\}$$
where the $E_i$ are finite nonempty sets (they are denoted $S_i$ in the paper  \cite{Gotfredsen-Zegers24}).
Denote by $N(v_i^s,v_j^t)$ the number of edges in $G$ from $v_i^s$ to $v_j^t$.
We have the following facts that are explained in detailed in \cite[Section 2]{Brzezinski-Szymanski18}:
\begin{itemize}
\item if $i>j$, then $N(v_i^s,v_j^t)=0$ for all $s,t$;
\item $N(v_i^s,v_i^t)=0$ if $s\neq t$;
\item $N(v_i^s,v_i^s)=1$ for each vertex $v_i^s$.
\end{itemize}
This implies that the graph $G$ satisfies the assumption of Proposition \ref{prop:spectrum}. Though, now there is $N:=\sum_{i=0}^n |S_i|$ vertices rather than $n$.

Consider the following example taken from \cite[Section 4]{Gotfredsen-Zegers24} where $N>n$.
Take the quantum lens space C$^*$-algebra $C^*(L_q(15, (2,2,5,2))$ which is a fixed point subspace of the C$^*$-algebra of a quantum 7-dimensional sphere.
This latter has for spectrum 4 circles.
Now, the associated graph $\Gamma$ has $1+1+5+1=8$ vertices (given by $|S_i|=gcd(15,5)=5$ if $i=3$ and $1$ otherwise) and thus the spectrum of the quantum lens space is parametrised by 8 circles rather than $4$.

\subsection{Even dimensional quantum sphere}
Fix $n\geq 1$ and denote by $S_q$ the quantum sphere of dimension $2n$ (with $0<q<1$).
As in the odd dimensional case Hong and Szyma\'{n}ski found an isomorphism between the C$^*$-algebra $C(S_q)$ and the Cuntz-Krieger C$^*$-algebra $\cO_{L}$ of a certain graph $L=L_{2n}$ described below, see \cite[Section 5.1]{Hong-Szymanski02}.
Recall that our convention demands that we take the graph of \cite{Hong-Szymanski02} but with the orientation of the edges reversed. This had no impact on the previous computations (since all the previous graphs were in fact isomorphic to their opposite) but now it will have its importance.

The graph $L$ has $n+2$ vertices denoted $\{1,\cdots,n,n+1,n+2\}$.
Here are the edges:
\begin{itemize}
\item for each $1\leq i\leq n$ there is one loop at the vertex $i$;\\
\item for each pair $(i,j)$ with $1\leq i<j\leq n$ there is one edge from $i$ to $j$;
\item for each $1\leq i\leq n$ there is one edge from $n+1$ to $i$ and one edge from $n+2$ to $i$.\\
\end{itemize}
The non-self-adjoint algebra $\cA:=\cA_L$ associated to this graph $L$ has for generators $n+2$ projections $a_1,\cdots,a_{n+2}$, for each edge from $i$ to $j$ an operator $a_{ji}$ satisfying $a_i a_{ji} a_j=a_{ji}$ and so that for each $j$ we have $\sum_i a_{ji}^* a_{ji}=a_j$.
More precisely, the isometric conditions translate into
\begin{itemize}
\item $a_{11}\oplus a_{1,n+1}\oplus a_{1,n+2}$ is a partial isometry with same domain than $a_1$;
\item $a_{21}\oplus a_{22} \oplus  a_{1,n+1}\oplus a_{1,n+2}$ is an isometry with same domain than $a_2$;
\item and more generally $\left(\oplus_{i=1}^j a_{ji}\right) \oplus a_{j,n+1}\oplus a_{j,n+2}$ is a partial isometry with same domain than $a_j$ for $1\leq j\leq n$. 
\end{itemize}

\begin{example}
Consider the case $n=3$. The graph has $5$ vertices and is as follows:
$$
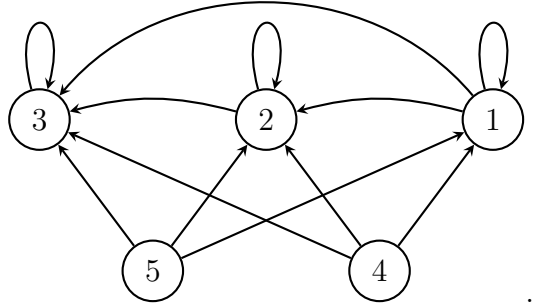
\begin{tikzpicture}[>=stealth,thick]

  \node[circle,draw,minimum size=8mm] (v3) at (0,2) {$3$};
  \node[circle,draw,minimum size=8mm] (v2) at (3,2) {$2$};
  \node[circle,draw,minimum size=8mm] (v1) at (6,2) {$1$};

  \node[circle,draw,minimum size=8mm] (v5) at (1.5,0) {$5$};
  \node[circle,draw,minimum size=8mm] (v4) at (4.5,0) {$4$};

  \draw[->] (v1) edge[loop above] node {} (v1);
  \draw[->] (v2) edge[loop above] node {} (v2);
  \draw[->] (v3) edge[loop above] node {} (v3);

  \draw[->] (v1) to[bend right=15] node[above] {} (v2);
  \draw[->] (v1) to[out=130,in=50] node[above] {} (v3);
  \draw[->] (v2) to[bend right=15] node[above] {} (v3);

  \draw[->] (v4) -- (v1);
  \draw[->] (v4) -- (v2);
  \draw[->] (v4) -- (v3);

  \draw[->] (v5) -- (v1);
  \draw[->] (v5) -- (v2);
  \draw[->] (v5) -- (v3);

\end{tikzpicture}.
$$
A module over $\cA$ is described by a Hilbert space $H$ that decomposes into $H=\oplus_{i=1}^5 H_i$ so that each arc from $i$ to $j$ defines a map 
$A_{ji}:H_j\to H_i$ and satisfying the axioms:
\begin{itemize}
\item $A_{11}\oplus A_{14}\oplus A_{15}:H_1\to H_1\oplus H_4\oplus H_5$ is an isometry;
\item $A_{21}\oplus A_{22}\oplus A_{24}\oplus A_{25}:H_2\to H_1\oplus H_2\oplus H_4\oplus H_5$ is an isometry;
\item $A_{31}\oplus A_{32}\oplus A_{33}\oplus A_{34}\oplus A_{35}:H_3\to H_1\oplus H_2\oplus H_3\oplus H_4\oplus H_5$ is an isometry.
\end{itemize}
A submodule is then a Hilbert subspace $H'$ of $H$ that decomposes into $\oplus_i H_i'$ with $H_i':=H'\cap H_i$ and satisfying that $A_{ji}(H_j')\subset H_i'$ for any edge-operator $A_{ji}.$
Note that \emph{any} Hilbert subspace of $H_4$ or of $H_5$ defines a submodule.
By taking $H'$ one-dimensional contained in $H_4$ or contained in $H_5$ we will then construct two irreducible modules. They will yield two irrreducible representations of the dimension $6$ quantum sphere. These two additional representations have no analogues in the representation theory of odd-dimensional quantum spheres.
\end{example}

Let us determine the finite dimensional irreducible modules.
Let $H$ be a finite dimensional irreducible module over the algebra $\cA$. 
We write $a_j\mapsto A_j$ and $a_{ji}\mapsto A_{ji}$ the action of $\cA$ on $H$.
Set $H_i$ the range (or domain) of the projection $A_i$ and recall that $A_{ji}$ can be interpreted as a map from $H_j$ to $H_i$.
There are no incoming edge to the vertex $n+1$. Therefore, $H_{n+1}\subset H$ is a submodule.
Moreover, note that since no operators of $\cA$ acts nontrivially on $H_{n+1}$ we have that any vector subspace of $H_{n+1}$ defines a submodule.
Irreducibility implies that $H_{n+1}$ is zero or one-dimensional.
In the latter case $H=H_{n+1}\simeq \C$ with $A_{ji}=0$ for all edges. Write $M_{n+1}$ this module structure.
Similarly, define $M_{n+2}$ to be the module where $H=H_{n+2}\simeq \C$ with $A_{ji}=0$ for all edges.
These are two irreducible modules (since they are of dimension one) and are inequivalent. 
Assume now that $H_{n+1}=H_{n+2}=\{0\}.$
This implies that all the operators $A_{j,n+1}$ and $A_{j,n+2}$ are equal to $0$ for all $1\leq j\leq n$.
Then $H_1\subset H$ is a submodule and thus $H_1=H$ or $\{0\}$ by irreducibility.
In the first case we have that $A_{11}:H_1\to H_1$ is an isometry and thus is a unitary since $H_1$ is finite dimensional by assumption.
Moreover, note that all other edge operators $A_{ji}$ are equal to $0$.
Hence, any invariant subspace of $H_1$ under $A_{11}$ defines a submodule.
Irreducibility implies that $H_1$ is one-dimensional and thus $A_{11}$ must be the multiplication by a scalar $z$ in the unit circle $S_1$.
Denote by $M_{1,z}$ this module structure, i.e.~$H=H_1\simeq \C$ with $A_{11}=z$ and $A_{ji}=0$ if $(j,i)\neq (1,1)$.
Assume now that $H_1=\{0\}$.
We may then apply a similar reasoning to $H_2$ and define a module $M_{2,z}$ for each $z\in S_1$.
Continuing inductively we deduce that a finite dimensional irreducible module is necessarily conjugate to $M_{n+1}$ or $M_{n+2}$ or $M_{j,z}$ for some $1\leq j\leq n$ and $z\in S_1$.
It is rather immediate that they are all indeed irreducible (since one-dimensional) and pairwise inequivalent.
Thoerem \ref{theo:Pi} implies that $\{ \Pi(M_{n+1}) , \Pi(M_{n+2}), \Pi(M_{j,z}):\ 1\leq j\leq n, z\in S_1\}$ provides a set of irreducible representations of $\cO$ that are pairwise inequivalent.
Moreover, if $K\simeq \Pi(H)$ with $\dim(H)<\infty$ so that $K$ is an irreducible representation of $\cO$, then $K\simeq \Pi(M)$ for $M=M_{n+i}$ or $M\simeq M_{j,z}$ for $i=1,2$ or $1\leq j\leq n,z\in S_1$.

Is the spectrum of $\cO$ larger?
Assume that $K$ is an irreducible representation of $\cO$ that is not isomorphic to $\Pi(M_{n+i})$ or $\Pi(M_{j,z})$ for $i=1,2$, $1\leq j\leq n$ and $z\in S_1$.
By essential surjectivity of the functor $\Pi$ there exists a module $H$ so that $K\simeq \Pi(H)$.
The module $H$ is necessarily indecomposable and, from the argument of above, does not contain any finite dimensional nonzero submodule.
As usual denote by $H_i$ the range of $A_i$ for $1\leq i\leq n+2$.
Note that $H_{n+1}\oplus H_{n+2}$ is a submodule as previously observed and any vector subspace of it defines a submodule.
Hence, $H_{n+1}\oplus H_{n+2}$ must be zero.
This implies that the edge operators $A_{j,n+1}$ and $A_{j,n+2}$ are equal to $0$.
We may then conclude that no such $K$ exists by applying a similar argument than in the odd-dimensional quantum sphere case. 

\begin{theorem}\label{theo:even-sphere}
Fix $n\geq 1$ and let $L=L_{2n}$ be the graph defined above with associated algebras $\cA,\cO$.
Identify $C(S_q)$ with $\cO$ where $S_q$ is the quantum sphere of dimension $2n$.
The set $$\{ \Pi(M_{n+1}) , \Pi(M_{n+2}), \Pi(M_{j,z}):\ 1\leq j\leq n, z\in S_1\}$$ constructed above which is parametrised by $n$ circles and two isolated points is a set of representatives of the spectrum of $C(S_q)$.
\end{theorem}

\end{document}